\documentstyle[12pt]{article}
\def\@themcountersep{}

\catcode`@=11
\renewcommand{\@begintheorem}[2]{\it \trivlist
      \item[\hskip \labelsep{\bf #1\ #2{.}}]}
\renewcommand{\@opargbegintheorem}[3]{\it \trivlist
      \item[\hskip \labelsep{\bf #1\ #2\ (#3) {:}}]}
\catcode`@=12
\newcommand{\proof}[1]{\hangindent=10pt \hangafter=0  \hskip 3 pt
               {\it Proof:\hskip 3pt} #1
              \hskip 3 pt \vrule height4pt  width3pt depth2pt
              \par\vskip 3pt plus 6pt
            }

\newtheorem{THEO}{Theorem}[section]
\newtheorem{ALGo}[THEO]{Algorithm}
\newtheorem{CONJ}[THEO]{Conjecture}
\newtheorem{COND}[THEO]{Condition}
\newtheorem{CORO}[THEO]{Corollary}
\newtheorem{DEFI}[THEO]{Definition}
\newtheorem{EXAMP}[THEO]{Example}
\newtheorem{FACT}[THEO]{Fact}
\newtheorem{HYPO}[THEO]{Hypothesis}
\newtheorem{LEMM}[THEO]{Lemma}
\newtheorem{PROB}[THEO]{Problem}
\newtheorem{PROP}[THEO]{Proposition}
\newtheorem{REMA}[THEO]{Remark}
\newcommand{\theo}{\begin{THEO}}
\newcommand{\algo}{\begin{ALGo} \rm}
\newcommand{\cond}{\begin{COND}}
\newcommand{\conj}{\begin{CONJ}}
\newcommand{\coro}{\begin{CORO}}
\newcommand{\defi}{\begin{DEFI} \rm}
\newcommand{\examp}{\begin{EXAMP} \rm}
\newcommand{\fact}{\begin{FACT}}
\newcommand{\hypo}{\begin{HYPO} \rm}
\newcommand{\lemm}{\begin{LEMM}}
\newcommand{\prob}{\begin{PROB} \rm}
\newcommand{\prop}{\begin{PROP}}
\newcommand{\rema}{\begin{REMA} \rm}
\newcommand{\etheo}{\end{THEO}}
\newcommand{\ealgo}{\end{ALGo}}
\newcommand{\econd}{\end{COND}}
\newcommand{\econj}{\end{CONJ}}
\newcommand{\ecoro}{\end{CORO}}
\newcommand{\edefi}{\end{DEFI}}
\newcommand{\eexamp}{\end{EXAMP}}
\newcommand{\efact}{\end{FACT}}
\newcommand{\ehypo}{\end{HYPO}}
\newcommand{\elemm}{\end{LEMM}}
\newcommand{\eprob}{\end{PROB}}
\newcommand{\eprop}{\end{PROP}}
\newcommand{\erema}{\end{REMA}}

\def\0{\mbox{\bf 0}}
\def\1{\mbox{\bf 1}}
\def\2{\mbox{\bf 2}}
\def\3{\mbox{\bf 3}}
\def\4{\mbox{\bf 4}}
\def\5{\mbox{\bf 5}}
\def\6{\mbox{\bf 6}}
\def\7{\mbox{\bf 7}}
\def\8{\mbox{\bf 8}}
\def\9{\mbox{\bf 9}}
\def\bxi{\mbox{\boldmath $\xi$}}
\def\etab{\mbox{\boldmath $\eta$}}

\def\cc{\mbox{\boldmath $c$}}
\def\d{\mbox{\boldmath $d$}}
\def\e{\mbox{\boldmath $e$}}

\def\q{\mbox{\boldmath $q$}}

\def\s{\mbox{\boldmath $s$}}

\def\u{\mbox{\boldmath $u$}}
\def\v{\mbox{\boldmath $v$}}
\def\w{\mbox{\boldmath $w$}}
\def\x{\mbox{\boldmath $x$}}
\def\y{\mbox{\boldmath $y$}}
\def\z{\mbox{\boldmath $z$}}

\def\I{\mbox{\boldmath $I$}}

\def\M{\mbox{\boldmath $M$}}

\def\Q{\mbox{\boldmath $Q$}}

\def\V{\mbox{\boldmath $V$}}

\def\Y{\mbox{\boldmath $Y$}}

\def\AC{\mbox{$\cal A$}}
\def\BC{\mbox{$\cal B$}}

\def\FC{\mbox{$\cal F$}}

\def\KC{\mbox{$\cal K$}}

\def\MC{\mbox{$\cal M$}}

\def\PC{\mbox{$\cal P$}}

\def\SC{\mbox{$\cal S$}}
\def\TC{\mbox{$\cal T$}}

\def\qf{\mbox{$q$\hspace{-0.15em}$f$}}

\def\NH{\mbox{$\widehat{\cal N}$}}
\def\DB{\mbox{$\overline{D}$}}

\begin{document}

\begin{center}

\begin{Large}
{\bf
Some
Fundamental Properties
of Successive Convex

\vspace{0.1cm}

Relaxation Methods
on LCP and Related Problems
}
\end{Large}

\bigskip

Masakazu Kojima \\
Department of Mathematical and Computing Sciences\\
Tokyo Institute of Technology \\
2-12-1 Oh-Okayama, Meguro-ku, Tokyo 152-8552, Japan \\
e-mail:{\it kojima@is.titech.ac.jp}

\bigskip

Levent Tun\c{c}el\footnote{Research of this author was supported
in part by a grant from NSERC and a PREA from
Ontario, Canada.
Some of this work was completed while this author
was a member of the Fields Institute, Toronto, Canada, during
September-December 1999.}\\
Department of Combinatorics and Optimization \\
Faculty of Mathematics \\
University of Waterloo \\
Waterloo, Ontario N2L 3G1, Canada \\
e-mail:{\it ltuncel@math.uwaterloo.ca}

\vspace{0.1cm}

October 8, 1999
\end{center}

\vspace{0.5cm}

\noindent
{\bf Abstract:}
General
Successive Convex Relaxation Methods (SRCMs)
can be used to compute the convex hull of
any compact set, in an Euclidean space,
described by a system of quadratic inequalities
and a compact convex set.
Linear Complementarity Problems (LCPs)
make an interesting and rich class of
structured nonconvex optimization problems.
In this paper, we study a few of the
specialized lift-and-project methods
and some of the possible ways of applying
the general SCRMs to LCPs and related
problems.

\vspace{0.1cm}



\noindent
{\bf Address for correspondence:}
Prof. Levent Tun\c{c}el,
Department of Combinatorics and Optimization,
Faculty of Mathematics,
University of Waterloo,
Waterloo, Ontario N2L 3G1, Canada.\\
e-mail:{\it ltuncel@math.uwaterloo.ca}

\newpage

\section{Introduction, SCRMs and Lov{\'a}sz-Schrijver procedures}

Since 1960's, complementarity problems attracted a very
significant attention in the theory as well as applications
of operations research. See, for instance, the
book on LCP \cite{CPS92}. In this paper, we consider various
complementarity problems in  the context of
Successive Convex Relaxation Methods (SCRMs)
proposed by the authors \cite{KT98a, KT98b}.
Since these methods can be used to compute
the convex hull of any compact subset of an
Euclidean space described by a system of quadratic inequalities
and a compact convex set, they can be used
to attack many complementarity problems
from several angles.

In the special case of 0-1 optimization problems
over convex sets, or more specially polytopes, there
are many Successive Convex Relaxation Methods (SCRMs)
based on lift-and-project techniques. We also
discuss some of the relationships of general SCRMs and these
more specialized algorithms in solving LCPs.

Let $F$ be a compact set in the $n$-dimensional Euclidean space
$R^n$. SCRMs take as input, a compact convex subset $C_0$ of $R^n$
and a set $\PC_F$ of quadratic functions which induce a description of $F$
such that
\[
     F = \{ \x \in C_0 : \qf(\x;\gamma,\q,\Q) \leq 0, \
             \qf(\cdot;\gamma,\q,\Q) \in \PC_F \}.
\]
Here we denote by $\qf(\cdot; \gamma, \q, \Q)$,
the quadratic function $(\gamma + 2\q^T\x + \x^T \Q \x)$.
Note that the variable $\x$ is irrelevant
outside a context and it will always be clear what the
variable vector is, from the context.

Let $\ell$ be an integer such that $1 < 2\ell \leq m$,
$\d \in R^m$, and let $A$ be a compact convex subset of $R^m$.
Consider the convex optimization problem with complementarity conditions:
\begin{equation}
\left.
\begin{array}{ll}
\mbox{maximize} & \d^T \u \\
\mbox{subject to} & \u \in A, \
                               0 \leq u_i, \ 0 \leq u_{i+\ell}, \
u_iu_{i+\ell} = 0, \
                                 \forall i \in \{1,2,\ldots,\ell\}.
\end{array}
\right\}
\label{ConvexLCP}
\end{equation}
First of all, it is clear that
LCP, with a known upper bound on a
solution of it, is a special case of (\ref{ConvexLCP})
(we can take $m = 2 \ell$ and $A$ as an affine subspace intersected
with a large enough ball).
Secondly, it is very elementary to formulate this
problem as a mixed 0-1 optimization
problem with convex constraints:
\begin{equation}
\left.
\begin{array}{lll}
\mbox{maximize} & \cc^T \v \\
\mbox{subject to} & \v \in C_0, \
                               v_i \in \{0, \ 1 \}, \
                                 \forall i \in \{m+1,m+2,\ldots,n\},
\end{array}
\right\}
\label{MixedIP}
\end{equation}
where
\begin{eqnarray*}
     C_0 & \equiv & \left\{\v =
     \left( \begin{array}{c} \u \\ v_{m+1} \\ \vdots \\ v_n\end{array} \right)
     \in R^{m+\ell} :
                          \begin{array}{lll}
                               \u \in A, \\
			0 \leq u_i \leq r v_{m+i}, \\
                               0 \leq u_{i+\ell} \leq r(1-v_{m+i}),\\
                               \forall i \in \{1,2,\ldots,\ell\}
                          \end{array}
                          \right\}, \
     \cc \equiv \left( \begin{array}{c} \d \\ \0 \end{array} \right) \in
R^{m+\ell}, \\
     n & \equiv & m + \ell, \
     r \geq \max_{i}\left\{ \max \{ u_i : \u \in A \}\right\}.
\end{eqnarray*}
In general, we allow $C_0$ to be an arbitrary compact convex
set in $R^n$.
There are various successive convex relaxation
methods that can be applied to such a problem.

We can represent the
feasible region $F \subset R^n$ of (\ref{MixedIP}) as
\[
      F = \{ \v \in C_0 : p(\v) \leq 0, \ \forall p(\cdot) \in \PC_F \},
\]
where $\PC_F$ denotes a set consisting of quadratic functions
\[
     (v_i^2 -v_i), \,\, (-v_i^2 + v_i), \ i \in \{m+1,m+2,\ldots,n\}
\]
on $R^n$.

In connection with the SCRMs and also the Lov{\'a}sz-Schrijver procedures
(see \cite{LS91}),
it seems convenient to introduce the following notation:
For every compact convex relaxation $C \subseteq C_0$ of $F$ and
every subset $D$ of $\DB \equiv \{ \d \in R^n : \|\d\| = 1\}$,
\begin{eqnarray*}
\PC^2(C,D)
& \equiv &
\{-(\d^T\v - \alpha(C_0,\d))(\bar{\d}^T\v - \alpha(C,\bar{\d})) :
\d \in D, \ \bar{\d} \in \DB \}, \\
\NH(C,D) & \equiv & \left\{ \v \in C_0 :
\begin{array}{l}
\exists \V \in \SC^n \ \mbox{such that } \\
\gamma + 2\q^T\v + \Q \bullet \V \leq 0,
\\ \forall \qf(\cdot;\gamma,\q,\Q) \in \PC_F \cup \PC^2(C,D)
\end{array}
\right\} \\
& & \mbox{(a Semi-Infinite LP relaxation of $F$)}, \\
\NH_+(C,D) & \equiv & \left\{ \v \in C_0 :
\begin{array}{l}
\exists \V \in \SC^n \ \mbox{such that }
\left(\begin{array}{ll} 1 & \v^T \\ \v & \V \end{array}\right) \in
\SC^{1+n}_+,
\\
\gamma + 2\q^T\v + \Q \bullet \V \leq 0,
\\ \forall \qf(\cdot;\gamma,\q,\Q) \in \PC_F \cup \PC^2(C,D)
\end{array}
\right\} \\
& & \mbox{(an SDP relaxation of $F$)},
\end{eqnarray*}
where $\alpha(C,\d) \equiv \max\{\d^T\v : \v \in C\}$
for every $\d \in \DB$.
Let $\SC^n$ and $\SC^{1+n}_+$ denote
the set of $n \times n$ symmetric matrices and the set of
$(1+n) \times (1+n)$ symmetric positive semidefinite matrices,
respectively. The corresponding variants of
Successive Semi-Infinite LP Relaxation Method
(SSILPRM) and
Successive SDP Relaxation Method (SSDPRM) can be
written as follows.

\algo (SSILPRM)
\label{SSILPRM}
\begin{description}
\item{Step 0: } Choose a $D_0 \subseteq \DB$. Let $k = 0$.
\item{Step 1: } If   $C_k = $ (the convex hull of $F$), then stop.
\item{Step 2: } Let $C_{k+1} = \NH(C_k,D_0)$.
\item{Step 3: } Let $k=k+1$, and go to Step 1.
\end{description}
\ealgo

\algo (SSDPRM)
\label{SSDPRM}
\begin{description}
\item{Steps 0, 1 and 3: } The same as the Steps 0, 1 and 3 of
Algorithm~\ref{SSILPRM}.
\item{Step 2: } Let $C_{k+1} = \NH_+(C_k,D_0)$.
\end{description}
\ealgo

To connect these algorithms to the
Lov{\'a}sz-Schrijver procedures,
we need to introduce some additional notation.
For every pair of closed convex cones
$\KC$ and $\TC$ in $R^{1+n}$, define
\begin{eqnarray*}
\MC(\KC,\TC) & \equiv &
\left\{
\Y = \left(\begin{array}{ll} \lambda & \lambda\v^T \\
           \lambda\v & \lambda\V \end{array}
\right) :
       \begin{array}{l}
          \lambda \geq 0, \ \v \in C_0, \ \V \in \SC^n, \\
            v_i = V_{ii}, \ i \in \{m+1,m+2,\ldots,n\}, \\
          \v^T\Y\w \geq 0, \ \forall \v \in \TC^*,
\ \forall \w \in \KC^*
       \end{array}
\right\}, \\
\MC_+(\KC,\TC) & \equiv &
\left\{
\Y = \left(\begin{array}{ll} \lambda & \lambda\v^T \\
                         \lambda\v & \lambda\V \end{array}
\right)  :
       \begin{array}{l}
          \lambda \geq 0, \ \v \in C_0,
\ \V \in \SC^n, \Y \in \SC^{1+n}_+\\
           v_i = V_{ii}, \ i \in \{m+1,m+2,\ldots,n\}, \\
          \v^T\Y\w \geq 0, \ \forall \v \in \TC^*,
\ \forall \w \in \KC^*
       \end{array}
\right\}.
\end{eqnarray*}
Let $\TC_0$ and $\KC_0$ be closed convex cones given by
\begin{eqnarray*}
      \TC_0^* & = & \mbox{c.cone} \left(\left\{\left(
\begin{array}{c} \alpha(C_0,\d)\\ -\d
\end{array}\right)  \in
R^{1+n}
: \d \in D_0 \right\}\right), \\
      \KC_0 & = & \left\{\left(\begin{array}{c} \lambda \\
\lambda\v \end{array} \right) \in R^{1+n} : \v \in C_0, \ \lambda
\geq 0 \right\}.
\end{eqnarray*}
(Note that $\TC_0$ itself is defined as the dual of $\TC_0^*$.)
If $C \subseteq C_0$ is a compact convex relaxation of $F$ and
\[
    \KC = \left\{\left(
\begin{array}{c} \lambda \\ \lambda\v \end{array} \right)
\in R^{1+m} : \v \in C, \ \lambda \geq 0 \right\},
\]
then
\begin{eqnarray*}
    \NH(C,D_0) & = & \left\{ \v \in R^n :
              \left(\begin{array}{ll} 1 & \v^T \\
                   \v & \V \end{array} \right) \in \MC(\KC,\TC_0)
             \right\}, \\
    \NH_+(C,D_0) & = & \left\{ \v \in R^n :
              \left(\begin{array}{ll} 1 & \v^T \\
                   \v & \V \end{array} \right) \in \MC_+(\KC,\TC_0)
             \right\}.
\end{eqnarray*}

Algorithms~\ref{SSILPRM} and~\ref{SSDPRM} specialized to
(\ref{MixedIP}) with
$\PC_F = \{  v_i^2 - v_i, \ -v_i^2 + v_i, \ i \in \{m+1,m+2,\ldots,n\} \}$
can be stated in the following forms,
which are essentially the Lov{\'a}sz-Schrijver procedures.

\bigskip

\noindent
{\bf Algorithm~\ref{SSILPRM}H}
 (Homogeneous form of Algorithm~\ref{SSILPRM}) \label{HSSILPRM}
\begin{description}
\item{Step 0: } Choose a $D_0 \subseteq \DB$.
Define $\TC_0$ and $\KC_0$ as above.
Let $k = 0$.
\item{Step 1: } If $\KC_k = \mbox{c.cone} \left(\left\{
\left(\begin{array}{c} 1 \\ \v \end{array} \right) : \v \in F
\right\}\right)$
then stop.
\item{Step 2: } Let
$\KC_{k+1} = \left\{ \Y\e_0 : \Y \in \MC(\KC_k,\TC_0) \right\}$.
\item{Step 3: } Let $k = k+1$, and go to Step 1.
\end{description}

\bigskip

\noindent
{\bf Algorithm~\ref{SSDPRM}H}
 (Homogeneous form of Algorithm~\ref{SSDPRM}) \label{HSSDPRM}
\begin{description}
\item{Steps 0, 1 and 3: }
The same as Steps 0, 1 and 3 of Algorithm~\ref{SSDPRM}H,
respectively.
\item{Step 2: } Let
$\KC_{k+1} = \left\{ \Y\e_0 : \Y \in \MC_+(\KC_k,\TC_0) \right\}$.
\end{description}

In this paper $\e_j$ denotes the $j$th unit vector and
$\e$ denotes the vector of all ones (the dimensions
of the vectors will be clear from the context).

\section{SCRMs applied to LCP
with an {\'a} priori bound}

Let $\M \in R^{\ell \times \ell}$,
$\q \in R^{\ell}$ be given. Consider the LCP
in the following form.
\begin{eqnarray*}
\mbox{(LCP) Find $\x$, $\s$ such that }
& & \M \x + \q \,\, = \,\, \s,\\
& & \x \,\, \geq \,\, \0, \,\, \s \,\, \geq \,\, \0,\\
& & x_i s_i \,\, = \,\, 0, \,\,
\forall \, i \in \{1,2, \ldots, \ell\}.
\end{eqnarray*}

Suppose we are given $\BC(\bxi, r) \equiv \left\{
\u \in R^{2\ell}:
\|\u - \bxi \| \,\, \leq \,\, r \right\}$,
an Euclidean ball containing a solution of the LCP.
(In the case of rational data $(\M, \q)$, we can take
$\BC$ centered at the origin with the radius bounded
above by a polynomial function of the ``bit size'' of
the data $(\M, \q)$.)
For the rest of this section, we assume that
the Euclidean ball with center $\bxi \equiv \0$
and the radius $r$ ($r$ is assumed given)
contains some solution of the LCP.

Under the boundedness
assumption above, it is particularly easy
to model any LCP as a 0-1 mixed integer programming
problem, since the only nonlinear constraints
of LCP can be expressed as
\[
x_i \,\, = \,\, 0, \,\,\,\, \mbox{ or }
\,\,\,\, s_i \,\, = \,\, 0, \,\,
\forall \, i \in \{1,2, \ldots, \ell\}.
\]
Balas' method \cite{Balas74} can be directly applied to such
formulations. We can also apply some variants of the
Lov{\'a}sz-Schrijver procedures \cite{LS91} to the mixed integer
programming feasibility problem:
\begin{eqnarray*}
\mbox{Find $\x$, $\s$ and $\z$ such that }
& & \M \x + \q \,\, = \,\, \s,\\
& & \0 \,\, \leq \,\, \x \,\, \leq \,\, r \z, \
\0 \,\, \leq \,\, \s \,\, \leq \,\,  r (\e - \z),\\
& & \z \in \{0,1\}^{\ell}.
\end{eqnarray*}

Note that we can eliminate the variable
vector $\s$ from the formulation
and apply the SSILPR and SSDPR Methods to the following
formulation:
\[
\begin{array}{l}
\0 \,\, \leq \,\, \M \x + \q \,\, \leq \,\,  r(\e - \z),\\
\0 \,\, \leq \,\, \z \,\, \leq \,\, \e,\
\0 \,\, \leq \,\, \x \,\, \leq \,\, r \z,\\
z_i^2 - z_i \,\, \leq \,\, 0, \,\,\,\, -z_i^2 + z_i \,\, \leq \,\, 0,
\,\,\,\, i \in \{1,2, \ldots, \ell\}.
\end{array}
\]
To apply the SCRMs, we can take
\begin{eqnarray*}
& & C_0 \,\, \equiv \,\, \left\{
\v= \left(\begin{array}{c} \x \\ \z
\end{array} \right) \in R^n: \begin{array}{l}
\0 \,\, \leq \,\, \M \x + \q \,\, \leq \,\,  r(\e - \z),\\
\0 \,\, \leq \,\, \z \,\, \leq \,\, \e,\
\0 \,\, \leq \,\, \x \,\, \leq \,\, r\z
\end{array}
\right\}, \\
& & m \equiv \ell,\ n \equiv 2\ell, \\
& & \PC_F \,\, \equiv \,\,
\left\{
(v_i^2 -v_i), \,\, (-v_i^2 + v_i), \ i \in \{m+1,m+2,\ldots,n\}
\right\}.
\end{eqnarray*}

Both algorithms, SSILPRM and SSDPRM presented in
Section 1,
terminate in at most $\ell$ steps.
This fact can be proved easily, using
the results of Balas \cite{Balas74},
Sherali and Adams \cite{SA90},
Lov{\'a}sz and Schrijver \cite{LS91},
or Kojima and Tun{\c c}el \cite{KT98a, KT98b}.
For computational experience on similar
algorithms for similar problems see
\cite{CP98}, \cite{SM96}.
In the next section, we give the details of a proof
of such a convergence result
when the methods are applied to a formulation
of Pardalos and Rosen \cite{PR88}.

\section{SCRMs applied to Pardalos-Rosen formulation of LCP}

We will illustrate the convergence proof on a formulation of (LCP) by
Pardalos and Rosen \cite{PR88}. They homogenize the vector $\q$
with a new continuous variable $\alpha$, then they maximize $\alpha$.
\begin{eqnarray*}
\left.
\begin{array}{lll}
(\mbox{MIP}_{\alpha}) & \mbox{maximize} & \alpha \\
& \mbox{subject to} &
\0 \,\, \leq \,\, \M \x + \q\alpha \,\, \leq \,\, \e-\z, \\
& & \0  \,\, \leq \,\, \x \,\, \leq \,\, \z,  \
 0  \,\, \leq \,\, \alpha \,\, \leq \,\, 1, \
 \z \,\, \in \,\, \{0,1\}^{\ell}.
\end{array}
\right\}
\end{eqnarray*}
Note that $\pmatrix{\bar{\alpha} \cr \bar{\x} \cr \bar{\z}}
\equiv \0$
is feasible in (MIP$_{\alpha}$) and, it is easy to see that
(MIP$_{\alpha}$) has an optimal solution
with $\alpha^* > 0$ iff
the (LCP) has a solution (or solutions) \cite{PR88}.
Moreover, if
$\pmatrix{\alpha^* \cr \x^* \cr \z^*}$ is an
optimal solution of (MIP$_{\alpha}$) with
$\alpha^* > 0$ then $\displaystyle \frac{\x^*}{\alpha^*}$ solves
the (LCP) \cite{PR88}.
One advantage of (MIP$_{\alpha}$) is that it does not
require the introduction of large, data dependent
constants (such as $r$ in the previous section)
or their a priori estimates.
Now,
we take
\begin{eqnarray*}
& & C_0 \equiv \left\{ \v =
\pmatrix{\alpha \cr \x \cr \z} \in R^{1+2\ell}:
\,\,
\begin{array}{l}
\0 \,\, \leq \,\, \M \x + \q \alpha \,\, \leq \,\,  \e - \z,\\
\0 \,\, \leq \,\, \x \,\, \leq \,\, \z,\
0 \,\, \leq \,\, \alpha \,\, \leq \,\, 1
\end{array}
\right\}, \\
& & m \equiv \ell +1, \ n \equiv 2 \ell +1, \\
& & \PC_F \,\, \equiv \,\,
\left\{
(v_i^2 -v_i), \,\, (-v_i^2 + v_i), \ i \in \{m+1,m+2,\ldots,n\}
\right\}.
\end{eqnarray*}

We have an analog of a
very elementary but also a key lemma (Lemma 1.3 of \cite{LS91})
of Lov{\'a}sz and Schrijver (and their proof technique
is adapted here).
In what follows, we refer to the vectors
in the space of $\KC_k$ by $\v$.
At the same time, we refer to different subvectors of
$\v$ by different names, such as $\x$, $\alpha$ etc., to
keep the correspondence of elements of $\v$
and the original formulation of $F$ clearer.
The proof of Lemma 1.3 of \cite{LS91} leads
to the following analogous result in our case.

\lemm \label{Lemma:LCP1}
Let $D_0 \supseteq \{ \pm \e_{m+1}, \pm \e_{m+2},
\ldots, \pm \e_{n} \}.$
Then the sequence of convex cones $\{ \KC_k: \,\, k \geq 0\}$
given by Algorithm 1.1H satisfies
\[
\KC_{k+1} \,\, \subseteq \,\,
\left(\KC_k \cap \left\{\v: x_i = 0 \right\}\right)
+ \left(\KC_k \cap \left\{\v: (\M\x +\q\alpha)_i = 0 \right\}\right),
\]
for every $i \in \{1,2, \ldots, \ell \},$
and for every $k \geq 0.$
\elemm

\proof{
Let $\w \equiv \pmatrix{1 \cr \bar{\alpha} \cr \bar{\x} \cr
\bar{\z}} \in \KC_{k+1}$. Fix
$j \in \{1,2, \ldots, \ell \}$ arbitrarily.
By the definition of $D_0$ and $\TC_0$, the unit vector
$\e_0$ is in $\TC_0$.
Hence, by the definition of $\MC(\KC_k, \TC_0),$
$\KC_{k+1} \subseteq \KC_k$ for every $k \geq 0$.
Therefore, $\w \in \KC_k$.
If $\bar{x}_j = 0$ or $(\M\x + \q\alpha)_j = 0$
then the statement of the lemma clearly holds.
So, without loss of generality, we assume
$\bar{x}_j>0$ and $(\M\x + \q\alpha)_j>0$.
Let $\Y \in \MC(\KC_k, \TC_0)$ such that $\w = \Y \e_0$.
By our choice of the cone $\TC_0$,
we conclude that $\Y\e_{n+j}$
and $\Y(\e_0 - \e_{n+j})$ are both in $\KC_k$.
Note that
\[
\w = \hat{\w} + \tilde{\w},
\]
where
$\hat{\w}\equiv \Y\e_{n+j}$
and
$\tilde{\w} \equiv\Y(\e_0 - \e_{n+j})$.
We will refer to the $\x$ and $\z$ parts of the vector $\hat{\w}$
by $\hat{\x}$, $\hat{\z}$ etc. (Similarly
for $\tilde{\w}.$)
First, since by the definition of
$\MC(\KC_k, \TC_0)$,
$v_{i} = V_{ii}$ for every $i \in \{m+1,
m+2, \ldots, n\}$, we have
$\tilde{z}_j = 0$ which implies
$\tilde{x}_j = 0$.
Therefore, $\tilde{\w}$ lies in the cone
$\left(\KC_0 \cap \left\{\v: x_j = 0 \right\}\right)$.
Second, since $\bar{x}_j > 0$, $\bar{z}_j$
must be positive. Therefore,
$(1/\bar{z}_j) \hat{\w} \in \KC_0$.
Since $v_{i} = V_{ii}$ for every $i \in \{m+1,
m+2, \ldots, n\}$,
$\hat{z}_j = \bar{z}_j$. So,
\[
\frac{1}{\bar{z}_j}\pmatrix{
\hat{\alpha} \cr
\hat{\x} \cr
\hat{\z}} \in C_k,
\]
with its $z_j$ entry equal to 1. Thus,
$(\M\hat{\x} + \q\hat{\alpha})_j = 0$.
Hence, $\hat{\w}$
is in the cone\\
$\left(\KC_k \cap \left\{\v: (\M\x +\q\alpha)_j = 0 \right\}\right)$.
Since the argument above is independent of the index $j$
the proof is complete.}

Note that the conclusion of the above lemma also applies
to the SSDPR Method since SSDPR Method yields at least as tight
relaxations as the SSILPR Method.

\theo \label{Thm:LCP1}
Both algorithms, Algorithm 1.1H and 1.2H
terminate in $\ell$ iterations when applied to the
formulation (MIP$_{\alpha}$) with our choice of $\PC_F$,
$C_0$ and $D_0$
above.
\etheo

\proof{
First note that
\[
\mbox{c.hull}(F) \subseteq \left\{\pmatrix{\alpha \cr
\x \cr \z} \in R^{n}:
\,\, \pmatrix{1 \cr \alpha \cr \x \cr \z} \in \KC_k \right\},
\,\,\,\, \forall \, k \geq 0.
\]
Next, let
$i,j \in \{1,2, \ldots, \ell\}$, $i \neq j$.
Since $\x \geq \0$ and
$\M\x + \q\alpha \geq \0$, for all $\v \in \KC_k$,
for every $k \geq 0$,
\begin{eqnarray*}
& & \left[\left(\KC_k \cap \left\{\v: x_i = 0 \right\}\right)
+ \left(\KC_k \cap \left\{\v:
(\M\x +\q\alpha)_i = 0 \right\}\right)\right]
\cap \left\{\v: x_j = 0\right\}\\
& = & \left(\KC_k \cap \left\{\v:
x_i = 0, \, x_j = 0 \right\}\right)
+ \left(\KC_k \cap \left\{\v:
x_j = 0, \, (\M\x +\q\alpha)_i = 0 \right\}\right).
\end{eqnarray*}
Similarly, for the intersection with
$\left\{\v: (\M\x +\q\alpha)_j = 0 \right\}.$
Now, we apply Lemma \ref{Lemma:LCP1}
repeatedly to conclude that
$\KC_{\ell}$ is the homogenization of the convex hull
of all solutions of the LCP that lie in the original relaxation
$C_0$.}

\section{SCRMs applied to the
smaller formulation of LCP
with explicit treatment of disjunctive constraints}

Now, we consider a formulation with fewer variables
and constraints.

\begin{eqnarray*}
\begin{array}{lll}
(\mbox{LCP}_{\alpha}) & \mbox{maximize} & \alpha \\
& \mbox{subject to} & \M \x + \q\alpha \,\, \geq \,\, \0, \,\, \x \,\, \geq
\,\, \0,
\,\, \alpha \,\, \geq \,\, 0,\\
& & \e^T(\M+\I)\x + (\e^T \q + 1) \alpha
\,\, \leq \,\, 1,  \\
& & x_i(\M\x+\q\alpha)_i = 0, \,\, i \in \{1,2, \ldots, \ell\}.
\end{array}
\end{eqnarray*}

It is easy to see that
$\pmatrix{\bar{\x} \cr \bar{\alpha}} \equiv \0$
is feasible in (LCP$_{\alpha}$),  and it is also
easy to observe that
(LCP$_{\alpha}$) has an optimal solution
with $\alpha^* > 0$ iff the (LCP) has a solution(or solutions).
Moreover, if
$\pmatrix{\x^* \cr \alpha^*}$ is an optimal solution of
(LCP$_{\alpha}$) with
$\alpha^* > 0$ then $\displaystyle \frac{\x^*}{\alpha^*}$ solves
the (LCP). Note that the inclusion in Lemma \ref{Lemma:LCP1}
can sometimes be strict for the SSILPR and SSDPR Methods.

We explicitly include the variable vector $\s$
in our discussion in this section,
for the sake of presentation.
Let
\[
C_0 \equiv \left\{\v
=\pmatrix{\x \cr \s \cr \alpha} \in R^{2\ell+1}: \,\,
\begin{array}{l}
\s \,\, = \,\, \M \x + \q\alpha \,\, \geq \,\,
\0, \, \, \x \,\, \geq \,\, \0, \,\, \alpha \,\, \geq \,\, 0, \\
\e^T(\M+\I)\x + (\e^T \q + 1) \alpha
\,\, \leq \,\, 1,\\

\end{array}
\right\}.
\]
In this section, we will describe another Successive
Convex Relaxation Method based on the ideas of
Balas \cite{Balas74}, Lov{\'a}sz and Schrijver \cite{LS91}.
This method will use only Linear Programming (LP) relaxations.
We describe
the method in the original space of $F$ and $C_0$.
Let $\FC(C_0)$ denote the set of facet
defining inequalities for $C_0$.
$\FC(C_0)$ is the input of the algorithm
which we introduce now.

\algo
\label{Alg:4.1}
{\bf Step 0.} $k \equiv 0.$\\
{\bf Step 1.} $\FC(C_{k+1}) \equiv \FC(C_k).$\\
{\bf Step 2.} For every inequality
\[
-\sum\limits_{i=1}^{\ell} \left(u_i x_i + u_{\ell+i} s_{i}\right)
- u_{2\ell+1} \alpha \,\, \leq \,\, u_0
\]
in $\FC(C_k)$ and every $j \in \{1,2, \ldots, \ell\}$
solve the LP problems
\begin{eqnarray*}
\begin{array}{lll}
(P_j) & \mbox{minimize} & \u^T \bxi^{(j)}\\
& \mbox{subject to} & \xi^{(j)}_j = 1,\
\xi^{(j)}_{\ell+j} = 0,\
\bxi^{(j)} \in \KC_k,
\end{array}
\end{eqnarray*}
and
\begin{eqnarray*}
\begin{array}{lll}
(P_{\ell+j}) & \mbox{minimize} & \u^T \bxi^{(\ell+j)}\\
& \mbox{subject to} & \xi^{(\ell+j)}_j = 0,\
 \xi^{(\ell+j)}_{\ell+j} = 1,\
 \bxi^{(\ell+j)} \in \KC_k.
 \end{array}
\end{eqnarray*}
If $(P_j)$ is infeasible then add the equation
$x_j = 0$ (or the inequality $x_j \leq 0$, since
the inequality $x_j \geq 0$ is already included)
to $\FC(C_{k+1})$.
If $(P_{\ell+j})$ is infeasible then add the equation
$s_j = 0$
to $\FC(C_{k+1})$.
Otherwise, let $(\bxi^{(j)})^*$ and $(\bxi^{(\ell+j)})^*$
denote the optimal solutions of $(P_j)$
and $(P_{\ell+j})$ respectively.
Define $y_j \equiv u_j - \u^T (\bxi^{(j)})^*$,
$y_{\ell+j} \equiv u_{\ell+j} - \u^T (\bxi^{(\ell+j)})^*$.
Add the inequality
\[
-\sum\limits_{i\neq j} \left(u_i x_i + u_{\ell+i} s_{i}\right)
-y_j x_j - y_{\ell+j}s_j -u_{2\ell+1} \alpha \,\, \leq \,\, u_0
\]
to $\FC(C_{k+1})$.\\
{\bf Step 3.} Let $k = k+1$, and go to Step 1.
\ealgo

Note that in iteration $k$, the algorithm
solves $(2\ell |\FC(C_k)|)$ LP problems.

\theo
\label{thm:4.1}
Let $C_k$, $k \in \{1,2, \ldots\}$ be the sequence
of convex relaxations generated by Algorithm \ref{Alg:4.1}.
Then $C_{\ell} = \mbox{c.hull}(F)$.
\etheo

\proof{
We think of
$\KC_k$  for all $k \geq 0$,
as a subset of $R^{1+(2\ell+1)}$, with the $0$th component
being the homogenizing variable, the next $\ell$ components
representing $\x$, the next $\ell$ components representing $\s$
and the last component representing $\alpha$.
Note that
\[
\KC_1 \,\, \subseteq \,\,
\left(\KC_0 \cap \left\{\v: x_j = 0 \right\}\right)
+ \left(\KC_0 \cap \left\{\v: s_j = 0 \right\}\right)
\]
iff
\begin{eqnarray}
\label{eqn:1}
\KC_1^* \,\, \supseteq \,\,
\left(\KC_0^* + \{-\e_j\}\right) \cap
\left(\KC_0^* + \{-\e_{\ell+j}\}\right).
\end{eqnarray}
(We used the fact that $\KC_0 \subseteq R_+^{1+(2\ell+1)}$.)
Therefore, if we ensure the latter inclusion,
then Theorem \ref{Thm:LCP1} applies and we can conclude
the
convergence of the method in $\ell$ iterations.
Recall that every vector $\u \in \KC_0^*$
represents a valid inequality
\[
-\sum\limits_{i=1}^{\ell} \left(u_i x_i + u_{\ell+i} s_{i}\right)
-u_{2\ell+1} \alpha \,\, \leq \,\, u_0
\]
for $C_0$.
To ensure the inclusion (\ref{eqn:1}),
it suffices to prove:
\[
\mbox{``For every } \u, \w \in \KC_0^*
\mbox{ such that }
u_i = w_i, \, \forall i \notin \{j, \ell+j\};
u_j \geq w_j, u_{\ell+j} \leq w_{\ell+j},
\]
\[\mbox{we have }
\y \in \KC_1^*,
\mbox{ where } y_i \equiv u_i, \forall i \neq j; \,\,
y_j \equiv w_j\mbox{.''}
\]
This is equivalent to proving the fact that if the two inequalities
\begin{eqnarray*}
& & -\sum\limits_{i=1}^{\ell} \left(u_i x_i
+ u_{\ell+i} s_{i}\right)
-u_{2\ell+1} \alpha \,\, \leq \,\, u_0, \
\mbox{ and} \\
& & -\sum\limits_{i\neq j} \left(u_i x_i + u_{\ell+i} s_{i}\right)
-w_j x_j - w_{\ell+j}s_j -u_{2\ell+1} \alpha \,\, \leq \,\, u_0
\end{eqnarray*}
are valid for $C_0$, then
\[
-\sum\limits_{i\neq j} \left(u_i x_i + u_{\ell+i} s_{i}\right)
-w_j x_j -u_{\ell+j}s_j -u_{2\ell+1} \alpha \,\, \leq \,\, u_0
\]
is valid for $C_1$.
To compute all such inequalities defining $C_1$,
we solve for every valid inequality
\[
-\sum\limits_{i=1}^{\ell} \left(u_i x_i + u_{\ell+i} s_{i}\right)
-u_{2\ell+1} \alpha \,\, \leq \,\, u_0
\]
for $C_0$ and every $j \in \{1,2, \ldots,\ell\}$,
the linear programming problems
\begin{eqnarray*}
\begin{array}{lll}
& \mbox{maximize} & \beta\\
& \mbox{subject to} &
\beta \e_j + \delta \e_{\ell+j} \,\, \preceq_{\KC_0^*} \,\, \u,
\end{array}
\end{eqnarray*}
and
\begin{eqnarray*}
\begin{array}{lll}
& \mbox{maximize} & \gamma\\
& \mbox{subject to} & \kappa \e_j + \gamma \e_{\ell+j} \,\,
\preceq_{\KC_0^*} \,\, \u.
\end{array}
\end{eqnarray*}
Here, $\preceq_{\KC_0^*}$ denotes the partial order induced by
the convex cone $\KC_0^*$ (that is, $\u^1 \preceq_{\KC_0^*} \u^2$
iff $(\u^2-\u^1) \in \KC_0^*$).
Note that both problems are always feasible.
Therefore, each of them either has an optimal solution or is unbounded.
If both LPs have optimal solutions, say $\beta^*$ and $\gamma^*$
then we set $w_j \equiv u_j - \beta^*$ and
$u_{\ell+j} \equiv u_{\ell+j} - \gamma^*$.
Since the above two problems are LPs, we can equivalently
solve their duals. Namely, we solve the LPs:
\begin{eqnarray*}
\begin{array}{lll}
(P_j) & \mbox{minimize} & \u^T \bxi^{(j)}\\
& \mbox{subject to} & \xi_j^{(j)} = 1,\
\xi^{(j)}_{\ell+j} = 0,\
\bxi^{(j)} \in \KC_0,
\end{array}
\end{eqnarray*}
and
\begin{eqnarray*}
\begin{array}{lll}
(P_{\ell+j}) & \mbox{minimize} & \u^T \bxi^{(\ell+j)}\\
& \mbox{subjcet to} & \xi_j^{(\ell+j)} = 0,\
 \xi_{\ell+j}^{(\ell+j)} = 1,\
 \bxi^{(\ell+j)} \in \KC_0.
 \end{array}
\end{eqnarray*}
These latter two linear programming problems are precisely
the ones used by Algorithm \ref{Alg:4.1}.
Notice that since their duals are either unbounded or
have optimal solutions, these LP problems either have
optimal solutions or are infeasible. When
$(P_j)$ is infeasible, the equality $x_j = 0$
is valid for $F$ and the algorithm adds this equality
to the describing inequalities of $C_k$.
Similarly, when $(P_{\ell+j})$ is infeasible,
$s_j = 0$ is valid for $F$ and the algorithm
behaves correctly in this instance. (In either
instance, the inclusion (\ref{eqn:1}) is obviously
satisfied for $j$.)
However, the proof is not yet complete; because, the arguments
so far ensure the inclusion
(\ref{eqn:1}) when the algorithm is ran for every
valid inequality of $C_0$. So, next we prove that
what the algorithm does (using only the facets of $C_0$) suffices.
To see this,
we need to prove that to derive the facets of $\KC_1$,
it suffices to start with a facet $\u$ of $\KC_0$
in the above procedure.
Suppose $\u, \w \in \KC_0^*$ satisfy the above
conditions but $\u$ is not facet inducing for
$\KC_0^*$. (We will prove that the valid inequality
derived from $\u$ and $\w$ is implied by
some other inequalities derived from some facets
$\u^1, \u^2, \ldots, \u^{\ell}$ of $\KC_0$.)
Since $\u$ is not facet inducing
for $\KC_0$, $\u$ is not an extreme ray of $\KC_0^*$.
Hence, there exist extreme rays
$\u^1, \u^2, \ldots, \u^{\ell}$ of $\KC_0^*$ such that
for some $\lambda_r > 0,$ $r \in \{1,2, \ldots,\ell\}$,
$\sum_{r=1}^{\ell} \lambda_r = 1$ the following
conditions are satisfied:
\begin{eqnarray*}
\u & = & \sum\limits_{r=1}^{\ell} \lambda_r \u^r,\
u_0^r  =   u_0, \,\, \forall r \in \{1,2, \ldots,\ell\}.
\end{eqnarray*}
Note that $\u^r$ is facet inducing for each $r$. Let $\bxi^r$
be the optimal solution of $(P_j)$
above for the objective function vector
$\u^r$. Let $\bxi^*$ be an optimal solution of $(P_j)$
when the objective function
vector is $\u$.
We claim that there exists $\tilde{\bxi} \in \KC_0$
such that
\begin{eqnarray*}
(\u^r)^T \tilde{\bxi} & = & (\u^r)^T \bxi^r, \,\,
\forall r \in \{1,2, \ldots,\ell\},\
\tilde{\xi}_j =  1, \
\tilde{\xi}_{\ell+j}  =  0, \
\tilde{\bxi}  \in \KC_0.
\end{eqnarray*}
(This claim follows from Farkas' Lemma, using the facts
that $\u^r \in \KC_0^*, \forall r$ and $\bxi^r \in \KC_0,
\forall r$.)
Thus, we have
\[
\sum\limits_{r=1}^{\ell}
\lambda_r (\u^r)^T \bxi^r \,\, = \,\, \u^T\tilde{\bxi}
\,\, \geq \,\, \u^T \bxi^*.
\]
Therefore, the inequality obtained from $\u$ is equivalent to
or dominated by a nonnegative combination of the
inequalities obtained from $\u^r$ which induce facets of $\KC_0$.
The proof is complete.}

We illustrated a derivation and convergence proof
for a successive relaxation method
(closely related to Balas' approach and
analogous to a suggestion of Lov{\'a}sz and Schrijver
\cite{LS91}) based on Lemma \ref{Lemma:LCP1}
and Theorem \ref{Thm:LCP1}.
Algorithm \ref{Alg:4.1} is an analog of a method
based on relaxations $N_0^k(\KC)$ from
\cite{LS91} (which is concerned with the case of
0-1 integer programming). For the relationship
of the methods of \cite{Balas74}
and \cite{LS91}, see Balas, Ceria and Cornuejols \cite{BCC93}.
(Balas' method \cite{Balas74}, in essence, corresponds to defining
\[
\KC_{k+1} \equiv \left(\KC_k \cap \left\{\v: x_{k+1} = 0 \right\}\right)
+ \left(\KC_0 \cap \left\{\v: (\M\x +\q\alpha)_{k+1} = 0 \right\}\right).)
\]

Let $C_k^{(4)}$, $k \geq 0$ denote
the projection of $C_k$ generated by Algorithm \ref{Alg:4.1}
onto the coordinates $\pmatrix{\x \cr \alpha}$.
Let $C_k^{(3)}$, $k \geq 0$ denote the projection
of $C_k$, generated by Algorithm 1.1, as used in Section 3,
onto the coordinates $\pmatrix{\x \cr \alpha}$.
Let $\KC_k^{(4)}$ denote the convex cone associated
with $C_k^{(4)}$.
{From} the proof of Theorem \ref{thm:4.1},
it is easy to see that
\[
\KC_{k+1}^{(4)} \,\, = \,\,
\bigcap_{i=1}^{\ell} \left[\left(\KC_k^{(4)} \cap
\left\{\v: x_i = 0 \right\}\right)
+ \left(\KC_k^{(4)} \cap \left\{\v: s_i = 0 \right\}\right)
\right].
\]
Therefore, the proofs of Theorems \ref{Thm:LCP1} and \ref{thm:4.1}
imply that
\[
\mbox{if } C_0^{(4)} \supseteq C_0^{(3)} \,\,\,\,
\mbox{ then } \,\,\,\, C_k^{(4)} \supseteq C_k^{(3)}
\,\, \mbox{ for all } k \geq 0.
\]
Thus, the SSILPR Method (Algorithm 1.1)
as applied in Section 3 to ($\mbox{MIP}_{\alpha}$)
converges at least as fast as Algorithm \ref{Alg:4.1}
applied to ($\mbox{LCP}_{\alpha}$).

\section{SCRMs applied to the
smaller formulation of LCP
with an implicit treatment of the disjunctive constraints}

We have already seen various ways of applying SCRMs to
LCP problems. Since the methods proposed in
\cite{KT98a, KT98b} only require a formulation of
the feasible solutions by quadratic inequalities, we are
also interested in applying the methods of \cite{KT98a, KT98b}
to the
following formulation:
\[
C_0 \equiv \left\{\pmatrix{\alpha \cr \x} \in R^{\ell+1}:
\,\,
\begin{array}{l}
\M \x + \q \alpha \,\, \geq \,\, \0, \,\, \x \,\, \geq \,\, \0,
\,\, \alpha \,\, \geq \,\, 0, \\
\e^T(\M+\I)\x + (\e^T \q + 1) \alpha
\,\, \leq \,\, 1
\end{array}
\right\},
\]
and
\[
\PC_F \equiv \left\{
x_i(\M\x + \q \alpha)_i \,\, \leq \,\, 0, \,\,\,\,
i \in \{1,2, \ldots, \ell\}
\right\}.
\]
The general theory of Kojima-Tun\c{c}el \cite{KT98a}
implies that their SSDPR and SSILPR Methods
converge.
It would be interesting to characterize the conditions
under which the Algorithms 3.1 and 3.2 of \cite{KT98b} converge
in at most $\ell$ iterations for the above
description of $\PC_F$ and $C_0$.
Also see \cite{KT99}, where the authors derived some necessary
and some sufficient conditions for the finite convergence of
SCRMs.

\section{A general linear complementarity problem}

Let $\AC: R^{\ell} \to R^{\ell}$,
a linear transformation,
$\q \in R^\ell$ and $\KC \subset R^{\ell}$ a pointed, closed
convex cone  with nonempty
interior, be given. Consider the Complementarity Problem
(CP):
\begin{eqnarray*}
\mbox{(CP) Find $\x$, $\s$ such that }
& & \AC(\x) + \q \,\, = \,\, \s,\\
& & \x \,\, \in \,\, \KC, \,\, \s \,\, \in \,\, \KC^*,\
\langle \x, \s \rangle \,\, = \,\, 0,
\end{eqnarray*}
where $\KC^*$ is the dual of $\KC$:
\[
\KC^* \,\, \equiv \,\, \{ \s \in R^{\ell} : \,\,
\langle \x, \s \rangle \,\, \geq \,\, 0,
\,\, \forall \, \x \in \KC \}.
\]
Since $\KC$ is a pointed, closed
convex cone with nonempty
interior, so is $\KC^*$.
Such problems
were studied recently, in the context of interior-point
methods \cite{SSK97}.
We pick $\etab \in \mbox{int}(\KC)$,
$\bar{\etab} \in \mbox{int}(\KC^*)$
and we can solve instead
the optimization problem
\begin{eqnarray*}
\begin{array}{lll}
(\mbox{CP}_{\alpha}) & \mbox{maximize} & \alpha \\
& \mbox{subject to} &
\x \,\, \in \,\, \KC,
\,\,\,\, \left[\AC(\x) + \q \alpha \right]
\,\, \in \,\, \KC^*, \,\,\,\, \alpha \,\, \geq \,\, 0,\\
& & \langle \bar{\etab}, \x \rangle + \langle \etab,
\AC(\x) + \q \alpha \rangle + \alpha \,\, \leq \,\, 1,\\
& & \langle \x, \AC(\x) + \alpha \q \rangle \,\, = \,\, 0.
\end{array}
\end{eqnarray*}

We choose
\[
C_0 \equiv \left\{\pmatrix{\alpha \cr \x} \in R^{\ell+1}:
\,\,
\begin{array}{l}
\x \,\, \in \,\, \KC, \,\,\,\,
\left[\AC(\x) +\q  \alpha \right]
\,\, \in \,\, \KC^*, \,\,\,\, \alpha \,\, \geq \,\, 0,\\
\langle \bar{\etab}, \x \rangle + \langle \etab,
\AC(\x) + \q \alpha \rangle + \alpha \,\, \leq \,\, 1
\end{array}
\right\}.
\]
Note that $C_0$ is always a compact convex set
(see the next theorem). We also pick
\[
\PC_F \,\, \equiv \,\,
\left\{
\begin{array}{c}
\langle \x, \AC(\x) + \alpha \q \rangle, \
-\langle \x, \AC(\x) + \alpha \q \rangle
\end{array}
\right\}.
\]

\theo
\label{thm:6.1}
\begin{itemize}
\item[(i)]
$C_0$ is a compact convex set.
\item[(ii)]
(CP$_{\alpha}$) has an optimal solution
with $\alpha^* > 0$ iff
(CP) has a solution (or solutions).
\item[(iii)]
If $\pmatrix{\alpha^* \cr \x^*}$ is an
optimal solution of (CP$_{\alpha}$) with
$\alpha^* > 0$ then the pair of vectors\\
$\displaystyle \left(\frac{\x^*}{\alpha^*},
\frac{1}{\alpha^*} \AC(\x^*) + \q \right)$ solves
(CP).
\end{itemize}
\etheo

\proof{
\begin{itemize}
\item[(i)] We only need to show that $C_0$ is bounded; because,
$C_0$ is a closed and convex subset of $R^{\ell+1}$ by definition.
Assume on the contrary that we can take an unbounded direction
$\displaystyle \pmatrix{\Delta\alpha \cr \Delta\x} \not= \0$ in $C_0$;
\begin{eqnarray*}
& & \pmatrix{\Delta\alpha \cr \Delta\x} \not= \0, \ \Delta\x \in \KC,
\Delta\alpha \geq 0, \
\left[\AC(\Delta \x) +\q  \Delta \alpha \right] \in \KC^*, \\
& & \langle \bar{\etab}, \Delta \x \rangle + \langle \etab,
\AC(\Delta\x) + \q \Delta \alpha \rangle + \Delta \alpha \,\, \leq \,\, 0.
\end{eqnarray*}
Since each term in the left hand side of the last inequality is
nonnegative, we have
\begin{eqnarray*}
\langle \bar{\etab}, \Delta \x \rangle = 0 \ \mbox{ and } \Delta \alpha = 0.
\end{eqnarray*}
Since $\bar{\etab} \in \mbox{int}(\KC^*)$ and $\Delta\x \in \KC$, the first
identity above implies that $\Delta\x = \0$. Thus,
we have a contradiction to
$\pmatrix{\Delta\alpha \cr \Delta\x} \not= \0$.
\item[(ii)]
Suppose (CP$_{\alpha}$) has an optimal solution
$\pmatrix{\alpha^* \cr \x^*}$
with $\alpha^* > 0$. Then
$\displaystyle \bar{\x} \equiv \frac{x^*}{\alpha^*} \in \KC$,
$\displaystyle \bar{\s} \equiv \frac{1}{\alpha^*}\AC(\x^*) + \q \in \KC^*.$
We have
\[
\langle \bar{\x}, \bar{\s} \rangle
\,\, = \,\, \langle \bar{\x}, \AC(\bar{\x}) + \q \rangle
\,\, = \,\, \frac{1}{(\alpha^*)^2}\langle \x^*, \AC(\x^*) + \alpha^*
\q\rangle = 0.
\]
Therefore, $(\bar{\x}, \bar{\s})$ solves (CP).
For the converse, let $(\bar{\x}, \bar{\s})$ be
a solution of (CP).
Let
\[
\zeta \equiv \langle \bar{\etab}, \bar{\x} \rangle
+ \langle \etab, \bar{\s} \rangle \,\, \geq \,\, 0, \
\alpha^* = \frac{1}{\zeta+1} \ \mbox{ and } \x^* = \frac{\bar{\x}}{\zeta+1}.
\]
Then $\pmatrix{\alpha^* \cr \x^*}$
is a feasible solution of (CP$_{\alpha}$).
But the feasible region of (CP$_{\alpha}$)
is compact and nonempty, its objective function is linear,
hence, (CP$_{\alpha}$) has optimal solution (or solutions).
Since we already showed a solution with positive objective value,
the optimum value is positive.
\item[(iii)]
This claim follows from the proof of (ii).
\end{itemize}
}

Theorem \ref{thm:6.1} shows that we can apply SCRMs
to (CP$_{\alpha}$) with the above $C_0$ and $\PC_F$
and solve the original, general problem (CP).

\end{document}